\documentclass[12pt]{article}
\usepackage{epsfig,amsmath,amssymb,stmaryrd,enumerate}
\usepackage{graphicx}

\newtheorem{theorem}{Theorem}
\newtheorem{col}{Corollary}
\newtheorem{rema}{Remark}

\begin{document}

\sloppy
\large

\title 
{One--point commuting difference operators of rank one
\thanks{The work is supported by RSF (grant 14-11-00441)}}
\author{Gulnara S. Mauleshova
\protect\footnote{
Novosibirsk state university; e--mail: mauleshova\_gs@mail.ru }, Andrey E. Mironov
\protect\footnote{Sobolev Institute of Mathematics, Novosibirsk, Russia;} \protect\footnote {Novosibirsk state university;
e--mail: mironov@math.nsc.ru}}
\date{}
\maketitle

\begin{abstract}
We consider one--point commuting difference operators of rank one. The coefficients of these operators depend on a functional parameter, shift operators being included only with positive degrees. We study these operators in the case of hyperelliptic spectral curve when the marked point coincides with the branch point. We construct examples of operators with polynomial and trigonometric coefficients. Moreover, difference operators with polynomial coefficients can be embedded in the differential ones with polynomial coefficients. This construction provides a new way of constructing commutative subalgebras in the first Weyl algebra.
\end{abstract}

If two difference operators
$$
L_k=\sum^{K^+}_{j=-K^-} u_j(n)T^j, \quad L_m=\sum^{M^+}_{j=-M^-} v_j(n)T^j, \quad n \in \mathbb Z,
$$
of orders $k$ and $m$, where $k=K^-+K^+,$ $m=M^-+M^+,$ $K^\pm,$ $M^\pm \geq 0,$ commute, then there exists a nonzero polynomial $R(z,w)$ such that $R(L_k,L_m)=0$ \cite{K3}. Polynomial $R$
defines the {\it spectral curve} $\Gamma=\{(z,w)\in \mathbb C^2 : R(z,w)=0\}.$ If $L_k\psi=z\psi, \quad L_m\psi=w\psi,$ then $P=(z,w) \in \Gamma.$ {\it Rank} $l$ of the pair $L_k,L_m$ is the dimension of the space of common eigenfunctions $l={\rm dim}\{\psi:L_k\psi=z\psi, \quad L_m\psi=w\psi\}$ for $P=(z,w) \in \Gamma$ in general position.
Any maximal commutative ring of difference operators is isomorphic to the ring of meromorphic functions on algebraic spectral curve with $s$ poles (see \cite{KN2}). Such operators are called {\it $s$--point operators}.
Eigenfunctions ({\it Baker -- Akhiezer functions}) of two--point rank one operators were found by I.M.~Krichever \cite{K3} (see also \cite{Mum}). Spectral data for the one--point operators of rank $l>1$ were received by I.M.~Krichever and S.P.~Novikov in \cite{KN2}. One--point rank two operators in the case of an elliptic spectral curve were found in \cite{KN2}. One--point rank two operators in the case of hyperelliptic spectral curve were studied in \cite{MM}.

Let us formulate our main results. We take the following spectral data
$$
S=\{\Gamma, \gamma_1,\ldots,\gamma_g, q, k^{-1}, P_n\},
$$
where $\Gamma$ is the Riemannian surface of genus $g$, $\gamma=\gamma_1+\cdots+\gamma_g$ is the non--special divisor on $\Gamma$, $q\in \Gamma$ is the marked point, $k^{-1}$ is the local parameter nearby $q$, $P_n\in\Gamma$ is the set of points, $n \in \mathbb Z$.

\begin{theorem}
There is a unique Baker -- Akhiezer function $\psi(n,P), \ n \in \mathbb{Z}, \ P\in\Gamma$ with the following properties.

\noindent{\bf 1}. The divisor of zeros and poles $\psi$ has the form
$$
\gamma_1(n)+\ldots+\gamma_g(n)+P_1+\ldots+P_n-\gamma_1-\ldots-\gamma_g-nq,
$$
in the case of $n \geq 0$ and
$$
\gamma_1(n)+\ldots+\gamma_g(n)-P_1-\ldots-P_n-\gamma_1-\ldots-\gamma_g-nq,
$$
in the case of $n < 0.$

\noindent{\bf 2}. In the neighborhood of $q$ the function $\psi$ has the expansion
$$
\psi=k^n+O(k^{n-1}).
$$
For any meromorphic functions $f(P)$ and $g(P)$ on $\Gamma$ with a single pole of orders $m$ and $s$ in $q$ with expansions
$$
f(P)=k^m+O(k^{m-1}), \quad g(P)=k^s+O(k^{s-1})
$$
there are unique difference operators
$$
L_m=T^m+u_{m-1}(n)T^{m-1}+\ldots+u_0(n),
$$
$$
L_s=T^s+v_{s-1}(n)T^{s-1}+\ldots+v_0(n)
$$
such that
$$
L_m\psi=f(P)\psi, \quad L_s\psi=g(P)\psi.
$$
The operators $L_m,L_s$ commute.
\end{theorem}

\begin{rema}
The spectral data with an additional set of points $P_n$ (analogously to our construction), were considered by I.M.~Krichever \cite{K1} in the case of two--dimensional discrete Schr\"{o}dinger operator.
\end{rema}

Notice that the divisor $\gamma_1(n)+\ldots+\gamma_g(n)$ is determined by the spectral data by the unique way. We also notice that in the special case where all points $P_n$ coincide, we obtain the two--point Krichever operators \cite{K3} of rank one.

Periodic two--point operators of rank one in which the shift operators have only negative degrees were discussed in a recent paper \cite{Krich}.

Consider the hyperelliptic spectral curve $\Gamma$ defined by the equation
\begin{equation}
\label{eq1}
w^2=F_g(z)=z^{2g+1}+c_{2g}z^{2g}+\ldots+c_0,
\end{equation}
with marked point $q=\infty.$ Let $\psi(n,P)$ be the corresponding Baker -- Akhiezer function. Then there exist commuting operators $L_2,L_{2g+1}$ such that
$$
L_2\psi=((T+U_n)^2+W_n)\psi=z\psi, \quad L_{2g+1}\psi=w\psi.
$$
\begin{theorem}
We have the equality
$$
L_2-z=(T+U_n+U_{n+1}+\chi(n,P))(T-\chi(n,P)),
$$
where
$$
\chi=\frac{\psi(n+1,P)}{\psi(n,P)}=\frac{S_n}{Q_n}+\frac{w}{Q_n},
$$
$$
S_n(z)=-U_nz^g+\delta_{g-1}(n)z^{g-1}+\ldots+\delta_0(n), \quad Q_n=-\frac{S_{n-1}+S_n}{U_{n-1}+U_n}.
$$
The functions $U_n, W_n, S_n$ satisfy the equation
\begin{eqnarray}
\label{eq2}
F_g(z)=S^2_n+Q_{n}Q_{n+1}(z-U^2_n-W_n)).
\end{eqnarray}
\end{theorem}

Equation (\ref{eq2}) can be linearized.

\begin{col}
The functions $S_n(z), U_n, W_n$ satisfy the equation
$$
(S_n-S_{n+1})(U_n+U_{n+1})-Q_n(z-U^2_n-W_n)+
$$
$$
Q_{n+2}(z-U^2_{n+1}-W_{n+1})=0.
$$
\end{col}

\begin{theorem}
In the case of an elliptic spectral curve $\Gamma$ given by the equation
$$
w^2=F_1(z)=z^3+c_2z^2+c_1z+c_0,
$$
operator $L_2$ of the form
$$
L_2=(T+U_n)^2+W_n,
$$
where
$$
U_n=-\frac{\sqrt{F_1(\gamma_n)}+\sqrt{F_1(\gamma_{n+1})}}{\gamma_n-\gamma_{n+1}},\quad W_n=-c_2-\gamma_n-\gamma_{n+1},
$$
$\gamma_n$ is an arbitrary function parameter, commutes with some operator $L_3.$
\end{theorem}

It can be shown that in the case of the hyperelliptic spectral curve and marked point $q=\infty$ operators $L_2,$ $L_{2g+1}$ can be obtained from one--point Krichever -- Novikov operators of rank two (see \cite{KN2}). We illustrate this in the case of $g=1.$ Under certain restrictions on the spectral data a one--point Krichever -- Novikov operator of rank two of order 4 with $g=1$ has the following form (this follows easily from \cite{MM})
$$
L_4=(T+U_n+V_nT^{-1})^2+W_n,
$$
where
$$
U_n=-\frac{\varepsilon_n+\varepsilon_{n+1}}{\gamma_n-\gamma_{n+1}},\quad W_n=-c_2-\gamma_n-\gamma_{n+1},
$$
$$
V_n=\frac{\varepsilon_n^2-F_1(\gamma_n)}{(\gamma_n-\gamma_{n-1})(\gamma_{n+1}-\gamma_{n})}.
$$
Operator $L_4$ commute with some operator $L_6=\sum_{j=-3}^3u_j(n)T^j.$
Coefficients of operators $L_4$ and $L_6$ are expressed through two functional parameters $\gamma_n, \varepsilon_n$. If one puts $\varepsilon_n=\sqrt{F_1(\gamma_n)},$ then one has the operators from theorem~3.

Theorem 2 enables us to construct explicit examples.

\begin{theorem}
The operator
$$
L_2=(T+r_1\cos(n))^2+\frac{1}{2}r_1^2\sec^2(g+\frac{1}{2})\sin(g)\sin(g+1)\cos(2n),
$$
$r_1\neq0$ commutes with operator $L_{2g+1}$ of order $2g+1$.
\end{theorem}

\begin{theorem}
The operator
$$
L_2=(T+\alpha_2n^2+\alpha_0)^2-g(g+1)\alpha_2^2n^2, \quad \alpha_2\neq0
$$
commutes with operator $L_{2g+1}$ of order $2g+1$.
\end{theorem}

\begin{rema}
One can directly check that if $g=1,\ldots,5$ then the operator
$$
L_2=(T+\alpha_2n^2+\alpha_1n+\alpha_0)^2-g(g+1)\alpha_2n(\alpha_2n+\alpha_1), \ \alpha_2\neq0
$$
commutes with $L_{2g+1}.$ Apparently this is true for all~$g$.
\end{rema}

Since
$$
[T,n]=T, \qquad [x,(-x\partial_x)]=x,
$$
the replacement $T \rightarrow x, \ n \rightarrow (-x\partial_x)$ in the operators
$$
L_2=(T+\alpha_2n^2+\alpha_1n+\alpha_0)^2-g(g+1)\alpha_2n(\alpha_2n+\alpha_1)
$$
and $L_{2g+1}$ yields a pair of commuting differential operators with polynomial coefficients and the operator $L_2$ corresponds to the operator
$$
(x+\alpha_2(x\partial_x)^2-\alpha_1(x\partial_x)+\alpha_0)^2-g(g+1)\alpha_2(x\partial_x)(\alpha_2(x\partial_x)-\alpha_1).
$$
So we obtain a commutative subalgebra in the first Weyl algebra $A_1=\mathbb{C}[x][\partial_x].$ The algebra $A_1$ has the following automorphisms $\varphi_j : A_1 \rightarrow A_1, \ j=1,2,3$
$$
\varphi_1(x)=\alpha x+\beta \partial_x, \ \varphi_1(\partial_x)=\gamma x+\delta \partial_x, \  \alpha,\beta,\gamma,\delta \in \mathbb{C}, \ \alpha \delta-\beta \gamma=1,
$$
$$
\varphi_2(x)=x+P_1(\partial_x), \quad \varphi_2(\partial_x)=\partial_x,
$$
$$
\varphi_3(x)=x, \quad \varphi_3(\partial_x)=\partial_x+P_2(x),
$$
where $P_1, P_2$ are arbitrary polynomials. J. Dixmier \cite{D} proved that the automorphism group $Aut(A_1)$ is generated by automorphisms of the form $\varphi_j.$ If applied $\varphi \in Aut(A_1)$ to $x, -x\partial_x \in A_1$, one obtains the elements $A=\varphi(x), \ B=\varphi(-x\partial_x),$ which also satisfy the equation
$$
[A,B]=A.
$$
The replacement $T \rightarrow A, \ n \rightarrow B$ in $L_2$ and $L_{2g+1}$ gives the commuting elements in $A_1.$ Thus the following important problem arises. To describe the solutions of the equation
$$
[A,B]=A, \quad A, B \in A_1
$$
up to the action of automorphisms $Aut(A_1)$. Each such solution allows one to construct commuting elements in $A_1$ via the commuting element in the ring of difference operators with polynomial coefficients $W_1=\mathbb{C}[n][T].$  As we told by P.S.~Kolesnikov, group of automorphisms $Aut(W_1)$ is generated by elements of the form
$
\varphi : W_1\rightarrow W_1,
$
$$
\varphi(T)=T, \qquad \varphi(n)=n+P(T),
$$
where $P$ is a polynomial.
Thus with the help of $Aut(W_1)$ and $Aut(A_1)$ one can obtains commuting differential operators via commuting difference ones with the same spectral curve. An interesting problem is to describe the commuting operators with polynomial coefficients with fixed spectral curve, which can be obtained from commuting difference operators via this procedure. This range of questions is associated with the Dixmier conjecture. This conjecture states that the $End(A_1)=Aut(A_1),$ or in the other words, if there is a solution of the string equation
$$
[A,B]=1, \ \ A,B \in A_1,
$$
then the operators $A$ and $B$ can be constructed from $\partial_x$ and $x$ with the help of some automorphism, i.e.
$$
A=\varphi(\partial_x), \quad B=\varphi(x), \quad \varphi \in Aut(A_1).
$$
The general Dixmier conjecture is stably equivalent to the Jacobian conjecture \cite{KK}. In a recent paper \cite{MZh} it was shown that the orbit space of the group action $Aut(A_1)$ on the set of commuting differential operators with polynomial coefficients with fixed spectral curve (\ref{eq1})
is always infinite if $g=1$ and for any $g$ there exists $F_g(z)$ with infinite number of orbits.
If one describes some class of commuting differential operators with polynomial coefficients (for example, derived from the difference operators) up to the action $Aut(A_1)$ with fixed spectral curve, it will give a chance to compare $Aut(A_1)$ and $End(A_1)$.

The authors are grateful to P.S.~Kolesnikov for useful discussions.


\end{document}